\documentclass[conference]{IEEEtran}
\usepackage{cite}
\usepackage[pdftex]{graphicx}
\usepackage{amsmath}
\usepackage{algorithmic}
\usepackage{array}
\usepackage{caption} 
\usepackage{amsmath}
\newtheorem{theorem}{Theorem}

\newtheorem{corollary}{Corollary}

\begin{document}
\IEEEoverridecommandlockouts

\title{Calibrating random number generator tests}
\author{\IEEEauthorblockN{Boris Ryabko}
\IEEEauthorblockA{Federal Research Center for Information and Computational Technologies
 of SB RAS\\ 
Novosibirsk state university
\\Novosibirsk, Russian Federation \\Email: boris@ryabko.net\\}
}

\maketitle

\begin{abstract}
Currently, statistical tests for random number generators (RNGs) are widely used in practice, and some of them are even included in information security standards.
But despite the popularity of RNGs, consistent tests are known only for stationary ergodic deviations of randomness (a test is consistent if it detects any deviations from a given class when the sample size goes to $ \infty $). However, 
 the model of a stationary ergodic source is too narrow for some RNGs, in particular, for generators based on physical effects. In this article, we propose computable consistent tests for some classes of deviations more general than stationary ergodic and describe some general properties of statistical tests. The proposed approach and the resulting test are based on the ideas and methods of  information theory.
\end{abstract}

\textbf{keywords:}   information theory, algorithmic information theory, randomness testing, statistical test,
random number generator,  Kolmogorov complexity, physical random number generators.

\section{Introduction and the statement of the problems}
The concept of randomness has attracted the attention of scientists for several millennia. In the 20th century, the concept of randomness and, in particular, of a random sequence began to be considered within the framework of algorithmic information theory (AIT) esteblished by A.~Kolmogorov,  
see \cite{Li:08,Calude:02,downey2010algorithmic,downey2006calibrating,merkle2006kolmogorov,vv,re,ve}.

Another line of research on randomness and randomness testing has been developed by researchers who develop practical random number generators  (RNGs) and tests for them, see 
\cite{le1,le2,NIST,j4}. It is worth noting that RNGs play an important role in cryptography, which is why some RNGs and randomness tests are included in the standards for data protection systems.  For example, the National Institute of Standards and Technology (NIST, USA) has proposed  the suite   (or ``battery'')  of 15 tests recommended as a standard \cite{NIST}.

Research in this area is conducted within the framework of the classical mathematical statistics
\cite{lem,ks}; specifically, it considers 
 a statistical test, a sample $\{ x = x_1x_2...\, x_n \} ,  n \ge 1,$ and 
  the following two statistical hypotheses:  
  \begin{equation}\label{H0}
   H_0 = 
\{ x  \emph{ 
is generated by the Bernoulli source with }
$$ $$
\emph{
probabilities (1/2, 1/2) \}      }  (B(1/2,1/2) )   
 \end{equation}
  and the alternative hypothesis $H_1$, which is a 
   set of probabilistic measures that does not contain $B(1/2,1/2)$.  
For example, the alternative hypothesis can be  
$H_1 = H^{st}_1$ which  is the set of all stationary ergodic measures, excluding $ B (1/2, 1/ 2) $.

 In mathematical statistics, the so-called consistent tests are highly desirable \cite{lem,ks};
 see also \cite{dr}. A formal definition will be given later, but here we note that 
 if $H_1$ is true and   such a test applies to the initial part $ x_1 ... x_n $ of $ x $, then it makes the correct decision from some $n$ on 
  when $ n \to \infty $ (that is, rejects $ H_0 $).

   Currently, among the hundreds of statistical tests applied to RNGs, consistent tests are known only for stationary ergodic alternatives (that is, $ H_1^{st}) $ and its subsets, see \cite {R2, r-isita}.  On the other hand, there may be real situations when the set of alternative hypotheses includes non-stationary sources. For example, RNGs popular in cryptography, which are  based on physical processes, may have some non-stationary deviations.
  This, in turn, indicates the need to search for tests for non-stationary deviations from randomness.

 So, the purpose of this article is to describe and analyze statistical tests that can be applied to an alternative hypothesis other than stationary deviations from randomness  with a focus on computable tests. 
We show that 
 general alternative hypotheses should be considered in the framework of AIT 
and describe some links between classical hypothesis testing and AIT.  
 Then   we  describe some consistent computable tests for the alternative hypotheses more general than  $ H_1^{st} $.
  In particular, it is done for $H_1$ which contains so-called $\sigma$-compressible sequences suggested in \cite{to} (see also \cite{rs}). This set of alternative hypotheses includes $H_1^{st}$ and, in a certain sense,  can be considered as its generalisation.
  The tests described are computable, but quite complex, and we show that some simple-to-implement tests can be effectively applied to broader alternatives than stationary deviations.
  
The rest of the article is organized as follows: the next part contains definitions and some auxiliary results. The third  part describes some general properties of    tests 
 in a framework of AIT. The  forth  part is devoted to the description and study of consistent computable tests for general alternative hypotheses, whereas the last part describes 
some simple-to-implement tests  that can be   applied beyond  stationary deviations.
The last part is a short conclusion.

\section{Definitions and auxilary results}\label{p1} 

\subsection*{The model.} 

 It is important to note that the general alternative hypothesis $H_1^{neg}  = \bar{H_0} 	$, that is, $H_1$ is negation of $H_0$, cannot be considered with the model 
 (\ref{H0}), because for any infinite sequence $x=x_1 x_2 ... $ there exists such a point measure  $\mu_x$, that $\mu_x(x) =1$. So, it is impossible to separate $H_0$ and $H_1^{neg} $ because any $x$  generated by $B(1/2,1/2)$ can be considered as generated by $\mu_x \in H_1^{neg} $.

To overcome this obstacle, one can attempt 
 to consider the following main hypothesis:
$H_0 = \{ x  \emph{ is  random  }\}   $
using the approach of 
  algorithmic information theory.

The notion of a random sequence plays a central role in algorithmic information theory, and there are currently several equivalent definitions within the AIT based on different  approaches \cite{re2}.   So, we propose to modify the main hypothesis as follows:  \begin{equation}\label{H0-2}
   H_0 = 
\{ x  \emph{ 
is  random }\},   
 \end{equation} 
and $H_1 $ is some subset of sequences which are not   random. 
Note that $\lambda (\{x; x \, is   \, random \} )= 1$, 
where $\lambda$ is the Lebesgue measure \cite{Li:08}.  In fact, it gives a possibility to apply  tests developed for (\ref{H0}) to  (\ref{H0-2}).

\subsection*{Statistical tests.}
Let there be the hypothesis $H_0$, 
 some alternative $H_1$,  
  $ T $ be a test and $\tau$ be a statistic  (that is, 
a function on   $ \{0,1 \}^n $ ) that is applied to a binary sequence $ x = x_1 ... x_n $.
(Here and below $\{0,1\}^n$ is the set of all $n$-bit binary words.)


By definition,  Type I error occurs   if
$H_0$ is true and  $H_0$ is  rejected and Type II error occurs   if
$H_1$ is true, but $H_0$ is accepted. 
The significance level $\alpha$, $\alpha \in (0,1)$, is the probability 
of the Type I error.  
Denote the critical region of the test $T$ for the  significance level $\alpha$ by $\bar{C}_T(\alpha,n)$  and let 
$C_T(\alpha,n)$ $ =  \{0,1\}^n \setminus \bar{C}_T(\alpha,n) .$ 
  Recall that, 
$H_0$ is rejected if and only if $x \in \bar{C}_T(\alpha)$ and, hence,
 \begin{equation}\label{alphasize}
|\bar{C}_T(\alpha,n)| \le 2^n \alpha \, .
 \end{equation}
 (Here and below $|X|$ is the number of elements in $X$ if $X$ is a set and the length of $X$ if $X$ is a word.)  
 Indeed, for both models,  
  if $H_0 $ is true, then all words in $\{0,1\}^n$ are  equiprobable (the probability of any word is $2^{-n}$). 

By definition, the test $T$ is consistent (for  $H_1)$, if for any $x \in H_1$ and any
significance level $\alpha \in (0,1)$, the hypothesis  
$H_0$ is rejected from some $n$ on 
 when $n \to \infty$. It  applies to both models (\ref{H0}) and (\ref{H0-2}), but  in the first case   it should be valid with probability one. 

Let us give a definition of the
 so-called p-value, which  plays an important role in the randomness testing \cite{NIST}. 
The  p-value  ($\pi_\tau (x)$) of  $\tau$ and  $x$  is defined by the equation
 \begin{equation}\label{pi-tau}
\pi_\tau (x) = \lambda \{ y: \tau(y) \ge \tau(x) \} = | \{ y: \tau(y) \ge \tau(x) \} | / 2^n \, .
 \end{equation}

 Informally, 
$ \pi_\tau (x)$ is the probability to meet a sequence  $y$ which is 
less random than the
observed $x$  when considering the null hypothesis. 
Clearly, 
$\bar{C}_T(\alpha) = \{x: \pi_\tau (x) \le \alpha \}$.

\subsection*{Batteries of tests.} 

Let us consider a situation where the randomness testing  is performed by conducting  a battery  of statistical  tests for randomness.  Suppose that the battery $\hat{T}$ contains a
finite or countable set of 
 tests   $T_1, T_2, ... $ and  $\alpha_i$ is the significance level of the $i-$th test, $i= 1,2,  ... $.  If the battery is applied in such a way that
 the hypothesis $H_0$ is rejected when at least one  test in the battery rejects it, then the significance level $\alpha$
of this battery satisfies the following inequality:
 \begin{equation}\label{a-g}
\alpha \le \sum_{i=1} \alpha_i \, .
\end{equation}
(Because $P(A+B) \le P(A)+P(B) $ for any events $A$ and $B$.) 

It will be convenient to formulate this inequality in a different way.
Suppose there  is  some $ \alpha  \in (0,1)$ and a sequence $ \omega$ of non-negative $ \omega_i$  such that $\sum_{i=1}^\infty
\omega_1 \le 1$. For example, we can define the following sequence $\omega^*$:
\begin{equation}\label{ome}
\omega^*_i =1 / (i (i+1)) \, \,\,\, i=1,2, ... \,\, .
\end{equation}
If  the significance level $T_i$ equals $\alpha \omega_i$, then the significance level of the battery $\hat{T}$ is not grater than $\alpha$.
(Indeed, from (\ref{a-g}) we obtain $\sum_{i=1} \alpha_i  = \sum_{i=1} (\alpha \omega_i)$ $ =\alpha \sum_{i=1} \omega_i \le \alpha $). 

If the p-value of the test $T_i$ is $\pi_{T_i}(\,)$, then
\begin{equation}\label{t-pi}
\pi_{\hat{T}}(\,) = \min_{i} \{\pi_{T_i} (\, ) / \omega_i \}. 
\end{equation} 
Note that this simple observation makes it possible to treat a  battery of tests as  a single test. 

\subsection*{Prefix-free codes,  consistent tests for stationary ergodic sources and universal codes.}
By definition, a set $U \subset  \{0,1\}^*$ is  prefix-free if for any different words $ u, v  \in U $ neither $ u $ is a prefix of $ v $, nor $v $ is a prefix of $  u $.
(Here and below  $ \{0,1\}^* $ and $ \{0,1 \}^\infty $ are sets of finite and one-sided infinite binary words, correspondingly.)
  The following Kraft inequality is valid for any prefix-free $U$:
 \begin{equation}\label{kraft}
\sum_{u \in U} 2^{-|u|} \, \le 1 \, .
\end{equation}
By definition, a map $\phi$ from   $\{0,1\}^*$ to  $\{0,1\}^*$  is a prefix-free code if for any different words $ u, v  \in \{0,1 \}^* $ $\phi(u) \neq \phi (v)$, 
and the set of all values $\phi(\,)$ is prefix-free \cite{co}. 

A code $\phi$ is universal if for any stationary ergodic $\nu$,  with probability one, 
 \begin{equation}\label{un-lim}
\lim_{n \to \infty } \frac{1}{n} | \phi( x_1 x_2 ... x_n) |  = h(\nu) \, ,
\end{equation}
where $ h(\nu) $ is the Shannon entropy of $\nu$ (see  \cite{co} for definition). Such codes exist;  see, for example, \cite{R2}.
Note that the purpose  
 of  codes is to ``compress'' sequences, i.e. to make the  length of the codeword 
$\phi( x_1 x_2 ... x_n) $ as small as possible. The   property (\ref{un-lim}) shows
that the universal codes are asymptotically optimal, because the Shannon entropy is 
a lower bound on the length of the compressed sequence (per letter), see \cite{co}.

It turns out that the following test is consistent for $H_1^{st}$ (\cite{R2,r-isita}):

{\it 
If  
 \begin{equation}\label{tau-fi}
 n -  |\phi (x_1 ... x_n)| \ge - \log_2 \alpha 
 \end{equation}
  then $H_0$ is rejected, 
otherwise it is accepted.
Here, as before, $\alpha$ is the significance level,   $\phi$ is  a universal code and
 $|\phi (x_1 ... x_n)|$  is the length 
of encoded (''compressed") sequence. }

It is worth noting that real data compressors (for example, zip)   can be used in this test, and such tests are quite effective \cite{rm}.
Moreover, such tests are asymptotically optimal in the following sense:  If we define the statistic 
$T_\phi (x_1 ... x_n) =  n -  |\phi (x_1 ... x_n)|$, then
 \begin{equation}\label{pi-univ}
 \lim_{n \to \infty} - \log \pi_{T_\phi} (x_1 ... x_n) = 1  - h(\mu) \, 
 \end{equation}
 and for any other statistic $\tau$  $\, \,\, \liminf_{ \to \infty} - \log \pi_{\tau} (x_1 ... x_n) \ge 1  - h(\mu) \,$,
 where $\mu $ is a stationary ergodic measure and $h(\mu)$ is its Shannon entropy, see \cite{r-isita}.

\subsection*{Some results of algorithmic information theory.}
First we  briefly introduce the notion of Kolmogorov complexity. A formal definition can be found, for example, in \cite{Li:08,Calude:02,downey2010algorithmic,downey2006calibrating}.  Informally, a prefix-free Kolmogorov complexity of a words $s$ is the length of the shortest program that outputs $s$. That is, for some universal prefix-free Turing machine $U$, we can define the Kolmogorov complexity $K_U(s) $ of a binary words as the length of the shortest program for $U$ that outputs $s$. (The Turing machine is prefix-free, if the set of $s$ on which $U$ halts forms a prefix-free code.)
There are such machines $U$ that $K_U(s)$ $\le K_U'(s)+const $ for every $s$ and every other universal Turing machine $U'$ (the constant may depend on  $U$ and $U'$, but not on $s$). Fix any such $U$ and define prefix-free Kolmogorov complexity $K(s)$ as $K_U(s)$. So, we can say that Kolmogorov complexity is defined up to an additive constant. 

The function $K(\, )$ is non-computable, but  there exists a computable function  $\hat{K} : N ,  \{0,1\}^* \to  \{1,2, ... \}$    such that for all integer $n$ and a binary word $ v$  : 
$\hat{K}(n, v) \le  \hat{K}(n + 1, v) $,  $  \lim_{n \to \infty } \hat{K}(n,v) = K(v)$
and for any $n$ the set $\hat{K}(n, v), v \in \{0,1\}^*$ is prefix-free
 \cite{Li:08}.

In what follows we will use the following definition of randomness which belongs to   Schnorr 
(and is equivalent to some other notions of randomness) 
 \cite{re2}: 
A sequence $x \in  \{0,1 \}^\infty$ is algorithmically random if there exists a constant $c$ such that, for all $n$,
\begin{equation}\label{shn}
K(x|_1^n) \ge n - c,
\end{equation}
 where $K$ is a prefix-free Kolmogorov complexity.

The following definitions are due to J.~Reimann \cite{re}.
 Let $x$ 
   be an infinite binary sequence (i.e. $x \in \{0,1\}^\infty$) and    $x|_1^n = x_1... x_n$.  Denote as  $0.x_1 ... x_n$   the real number whose binary representation is the real number $0.x_1 ... x_n$ from $[0,1)$. 
The lower pointwise dimension for a measure $\mu$ and $x = x_1x_2 ... \in \{0,1 \}^\infty$  is defined by
\begin{equation}\label{t-pi}
\delta_\mu(x) = \liminf_{n \to \infty} \frac{-1}{\, \, n}
 \log \mu([0.x_1 ... x_n ,
0.x_1 ... x_n + 2^{-n}) )
\end{equation}

It is shown in AIT that 
 there exist a so-called universal semimeasure $M$ which is enumerable, and for any computable measure $\mu$ there exists such a constant $C$ that  the following inequality holds \cite{re}:
\begin{equation}\label{Kmu}
| \, \, -\log M(u)  + \log \mu (u) | \, \le C .
\end{equation}

\subsection*{$\sigma$-compressible sequences and tests for stationary ergodic sources.}
The set of  $\sigma$-compressible  sequances is defind in \cite{to}.  By definition,  for $\sigma \in [0,1]$,  $x \in \{0,1 \}^\infty$ is $\sigma$-compressible if 
\begin{equation}\label{si-co}
K(x|_1^n) \le \sigma n +o(n) \, ,
\end{equation}
where $n$ grows. Let $\Sigma$  denote the set of  
 all sequences which are $\sigma$-compressible for a $\sigma \in [0,1]$.

It can be seen that for any stationary ergodic measure $\mu$  the set of  $\sigma$-compressible sequences contains the set of $\mu$-typical sequences,
if $h(\mu) = \sigma$. So,   
$\mu( \Sigma) = 1$ \cite{rs} and, in a    sense, the notion of $\sigma$-compressibility is more general than that of a  stationary ergodic source.  

\section{Links between classical statistics and AIT }
It turns out that some central concepts of mathematical statistics are deeply related to the results of the AIT.
They are summarized in the following theorem.
 \begin{theorem}\label{t2} Let $H_0 = \{ x  \emph{ is  random  }\}   $
and  $\mu$ be a measure for which the hypothesis $H_1$ in  (\ref{H0})  is true (i.e. $\mu \neq B(1/2, 1/2)$).

i) Then, with $\mu$-probability 1, for $x= x_1 x_2 ...$    
$$ \limsup_{n \to \infty } \frac{- \log \pi_\mu(x|_1^n)  }{n} \ge 1 - \delta_\mu(x) \, ,$$
where $\pi_\mu(x|_1^n)  $ is the p-value for the statistic $\tau(x) = \mu(x)$.

ii) For a universal semicomputable measure $M$ and a computable $\mu$, with $\mu$-probability 1, 
$$\limsup_{n \to \infty } \frac{- \log \pi_M(x|_1^n)  }{n} \ge 1 - \delta_\mu(x) ,$$ 

\end{theorem} 
  (In a sense, i) corresponds to the case of the known hypothesis 
   $ H_1 $, and  ii) to the case of an unknown $H_1$.)
  
  {\it Proof of the theorem}.  By definition (\ref{pi-tau}),
  $$
  \pi(\mu(x|_1^n)) = | \{y_1 ... y_n: \mu(y_1 ... y_n) \ge \mu(x|_1^n) \}|/ 2^n \, .
  $$
Clearly, $| \{y_1 ... y_n: \mu(y_1 ... y_n) \ge \mu(x|_1^n) \}| \le 1/ \mu(x|_1^n)$.
From the two latter inequalities we obtain 
$\pi(\mu(x|_1^n)) \le 1/ (\mu(x|_1^n) 2^n) .$
So,  $ \log \pi( \mu(x|_1^n)) \le  - n -  \log \mu(x|_1^n)$.
Hence, 
$$
\limsup_{n \to \infty } \frac{- \log \pi_\mu(x|_1^n)  }{n} \ge 1 -
\liminf_{n \to \infty } \frac{- \log \mu(x|_1^n)  }{n} \, .
$$
Taking into account the definition of the lower pointwise dimension (\ref{t-pi}) we obtain the first statement of the theorem.

The second statement of the theorem immediately follows from (\ref{Kmu}).

\section{Computable tests for general alternative hypotheses }

In this part we will consider the general hypothesis $H_0$ (\ref{H0-2}). 
It is known in AIT that a computable test for the most general alternative hypothesis 
$H_1^{neg} =$  $\bar{H}_0$ does not exist. So, we consider two alternative hypotheses for which computable tests exist.
Informally,
suppose that $\omega $ is  such a sequence of non-zero $\omega_i, i=1, 2 , ...  $
 that $\sum_{i=1}^\infty \omega_i \le 1$.  The first alternative hypothesis    is as follows:
   $H_1^\omega  = \{x : K(x|^n_1 > n - \log (1/\omega_n) \}$. 
In particular,  if $ \omega_n = \it{const} /n^{1+\epsilon}$ for some $\epsilon >0$, then 
$H_1^\omega  = \{x : K(x|^n_1 > n - (1+\epsilon )\log n + O(1) \}$, if 
$ \omega_n =  \it{const}/(n \log n^{1+\epsilon})$, then 
$H_1^\omega  = \{x : K(x|^n_1 > n - ( \log n + (1+\epsilon )\log \log (n+1))  +O(1) \}$, etc. 

It is interesting to compare these equations with the Schnorr definition of a random sequence (\ref{shn}).
In a sense, the $ O (\log n) $ distance is the price to pay for the computability of tests.

 Now let's continue the description of the test.  Let $f(n), n \in N,$ be a monotonically increasing  computable function such that $f(n) \to \infty$.
Define the following  statistic  and test:
\begin{equation}\label{test}
\tau_K (x_1... x_n) = \min_{m = 1, ..., n} (m - \hat{K}(f(n), x|_1^m) - \log(1/\omega_m) ) ,
$$ $$
H_0 \, \,  is \, \,  rejected \, \,  if  \, \, \tau_K (x|_1^n) \ge \log (1/\alpha)  ,
\end{equation}
 where $\alpha $ is the significance level.
 Informally, this  is a battery of $n$ tests $m - \hat{K}(f(n), x|_1^m) - \log(1/\omega_m)$, $m = 1, ... , n$, where $\alpha$ is "shared" among the tests. When $n$ grows, $f(n)$ 
 goes to $\infty$ and, hence, $ \hat{K}(f(n), x|_1^m)  $  $\to K(x|_1^m)$. The following theorem is based on this  observation.
 \begin{theorem}\label{t1}
Let  there be 
a sequence $x = x_1 x_2 ... $ such that 
\begin{equation}\label{cond}  \liminf_{n \to \infty} ( n - K(x|_1^n) - \log (1/\omega_n)  )= \infty.  
 \end{equation}
 Then the   test (\ref{test}) is consistent for $x$. 
\end{theorem}
\begin{corollary}\label{t1}
The test  (\ref{test}) is consistent for any $\sigma$- compressible sequence $x$ if 
$\sigma < 1$.
\end{corollary}
{\it Proof of the theorem 2}. 
First we show that the significance level is not grater than $\alpha$. 
To do that, 
consider the statistic $\tau^*(m,n)$ and the test  $T^*(m,n)$ defined as follows:  
\begin{equation}\label{test2} 
\tau^*(m,n) = m - \hat{K}(f(n), x|_1^m) - \log(1/\omega_m) \,  ,
$$ $$
H_0 \, \,  is \, \,  rejected \, \,  if  \, \, \tau^*(m,n) \ge \log (1/\alpha)  .
\end{equation}
 Suppose that $H_0$ is true and let $ 1\le m \le n$.  From the definition of the test 
 (\ref{test2}) we can see that the probability that $H_0$ is rejected is 
 $|W|/2^m$, where  
$W = \{w: w \in \{0,1\}^m \, and \, $ $m- \hat{K}(f(n), w ) - \log(1/\omega_m) $ $\ge 
\log(1/ \alpha) \}$. Clearly, $W= \{ w:  w \in \{0,1\}^m \, and \, $ $\hat{K}(f(n), w ) \le m - \log (1/ (\alpha \omega_m)) \}$. For any $n$ the set
$\hat{K}(f(n), v), v \in \{0,1\}^m$ is prefix-free.
Taking into account the definition of $W$,   from the Kraft inequality we obtain 
$$
1 \ge \sum_{w \in W} 2^{-|w|} \ge |W| 2^{- (m - \log (1/ (\alpha \omega_m)) )} \, .
$$
Hence, $ |W| \le 2^{m + \log (1/ (\alpha \omega_m))}$.
So, the probability that $H_0$ is rejected (if it is true) equals $|W|/2^m \le \alpha \omega_m.$

From definitions  (\ref{test}) and  (\ref{test2}) we can see that the first of them is the battery of the test 
and from (\ref{a-g}) we can see that  the probability that $H_0$ is rejected (if it is true) is not grater than $\sum_{i=1}^n  \alpha \omega_i \le  \alpha$.

We can see from (\ref{cond}) that there exists  such a increasing sequence $m_1, m_2,  ... $  that $\lim_{k \to \infty} ( m_k - K(x|_1^{m_k}) + \log (1/\omega_{m_k})  )= \infty $. Hence, for any $\alpha \in (0,1)$ there exists such $m^*$ that 
\begin{equation}\label{alph} ( m^*- K(x|_1^{m^*}) + \log (1/\omega_{m^*})  ) > \log (1/\alpha ) \, .
\end{equation}  
By definition, 
$  \lim_{l \to \infty } \hat{K}(l,x|_1^{m^*}) = K(x|_1^{m^*})$ and, hence,
$  \lim_{n \to \infty } \hat{K}(f(n),x|_1^{m^*}) = K(x|_1^{m^*})$. 
Taking into account (\ref{alph})  from the last equation and the definition of the statistic
(\ref{test}) we can see that there exists such $n^*$ that  $\tau_K(x|_1^n) >
\log(1/\alpha)$ if $n > n^*$.

{\it Proof of the corollary }. 

Taking into account the definition (\ref{si-co}), we can see that  for any $\sigma$-compressible $x$ and $\sigma < 1$ 
$$  
n - K(x|_1^n) \ge (1-\sigma) n + o(1) \,  .
$$
So, for any $\alpha \in (0,1)$, there exists such an integer $m$ that 
$$
K(x|_1^m) > \log( 1/ (\omega_m  \alpha ) ) \, .
$$
From this and the definition of $\hat{K}$ and $f$ we can see that 
$$ \ lim_{ n \to \infty}  \hat{K}(f(n), x|^1_m)  > \log( 1/ (\omega_m  \alpha ) ) \, ,
$$
and, hence, the test (\ref{test}) rejects $H_0$.

\section{Simple-to-implement tests consistent to some 
 non-stationary deviations from randomness}
In this part we will describe some easy-to-implement tests that can be applied to some of the alternatives outside of $ H_1^{st} $. 

First, we will briefly describe the LZ77 code based on the definition in \cite{gram}.
Suppose, there is a binary string $\sigma^*$ which is encoded using the code LZ77. 
This string  is represented by a list of pairs
$(p_1; l_1) ...  (p_s; l_s)$. Each pair $(p_i; l_i)$ represents a string,
and the concatenation of these strings is $\sigma^*$. In particular,
if $p_i = 0$, then the pair represents the string $l_i$, which is
a single terminal. If $p_i  \neq 0$, then the pair represents a
portion of the prefix of $\sigma^*$ that is represented by the preceding
$i - 1$ pairs; namely, the $l_i$ terminals beginning at
position $p_i$ in $\sigma^*$; see \cite[part 3.1]{gram}.
The length of the codeword depends on the encoding of the sub-words $p_i, l_i$
which are integers. For this purpose we will use a prefix code $C$ for integers, for which for any
integer $m$  
\begin{equation}\label{C}
   |C(m)| = \log m +2 \log \log (m +1) +O(1)  \,\,.
 \end{equation}
 (Such  codes are known in Information Theory, see, for example,  \cite[part 7.2]{co}).
 Note that $C$ is the prefix code and, hence,  for any $r \ge 1$ the codeword $C(p_1)C(l_1) ... C(p_r) C(l_r)$ can be decoded to $(p_1;l_1) ... (p_r; l_r)$.
There  is the following upper bound for the length of the LZ77 code \cite{co,gram}: for any
word $w_1 w_2 .... w_m$  
\begin{equation}\label{leng2}
   | code_{LZ}(w_1 w_2 ... w_m)| \le m \, (1+o(1) ),
 \end{equation}
 if $m \to \infty$.

We recall the definition of a typical sequence for the stationary ergodic measure $ \mu $.
The sequence $ x  \in \{0,1\}^\infty$ is $ \mu $ -typical if the limiting frequency of occurrence of any word $ u \in \{0,1 \} $ in $ x $ is $ \mu (u) $. 
Besides, 
 the  $\mu$  measure of the set of all $\mu$-typical sequences equals 1. 

We will now describe such   sequences that, on the one hand,   are not  typical for any stationary ergodic measure and, on the other hand,  are not   random and will be rejected by the suggested test. Thus, the proposed model makes it possible to detect non-random sequences whose deviantions from randomness are not stationary, 
i.e., which are outside $ H_1^{st} $.  Those sequences  are recognized by LZ77 and some grammatical codes as non-random, if those codes are used in the test 
(\ref{tau-fi}). 
To do this, we take any binary    incompressible 
sequence $ x = x_1 x_2 ... $  
(that is, for which (\ref{shn}) is valid) 
and define a new sequence $ y (x) = y_1 y_2. .. $ as follows. Let for   $k =0,1, 2, ... $ 
 \begin{equation}\label{Y}
  u_k = x_{2^{2^k}-1} x_{2^{2^k}} x_{2^{2^k}+1 }... x_{2^{2^{k+1}}-2}
$$ $$
y(x)= u_0 u_0 u_1 u_1  u_2  u_2 u_3 u_3  ...
\end{equation}
(For example, $u_0 = x_1 x_2 $, $u_1 = x_3 x_4 ... x_{14}$, $u_2 = x_{15} ... x_{254}$, 
$y=  x_1 x_2  $ $x_1 x_2   x_3 x_4 ... x_{14} $ $ x_3 x_4 ... x_{14}$ 
$x_{15} ... x_{254}$ $x_{15} ... x_{254}$ $...$. 

The idea behind this design is quite clear. Firstly, it is obvious that the word $y$ 
cannot be typical for a stationary ergodic source and,  
secondly, when $u_0 u_0 u_1 u_1 ...  u_k u_k$ is encoded 
 the second subword $ u_k$ will be encoded by a very short word (about $ O (\log |u_k|))$, since it coincides with the previous word  $u_k$. So, for large $k$ the length of 
 the encoded word $u_0 u_0 u_1 u_1  ...  u_k u_k$ will be about 
 $|u_0 u_0 u_1 u_1 ...  u_{k} u_k| \, (1+ o(1) )/\, 2$. In turn, 
 taking into account the test (\ref{tau-fi}) we can see that  $H_0$ will be rejected for any $\alpha \in (0,1)$ when $k$ grows. The following is a more formal statement based on this simple consideration: 

{\bf Claim.}  
For any incompressible 
sequence $x$  the sequence $y(x)$ is not typical for any stationary ergodic measure and 
the test (\ref{tau-fi}) along with 
the LZ77  code rejects $x$ for any significance level $\alpha \in (0,1)$.

\section{Conclusion}  

In this article, we have shown that such an applied field as RNG testing requires new models in addition to the probabilistic approach and, in particular, in addition to  the stationary ergodic model. Naturally, suitable models can be found within algorithmic information theory, which, in a sense, was developed for the study of randomness and, hence, RNGs.
 It turned out that the gap between computability and non-computability of consistent tests is about $log $ distance in Kolmogorov complexity.

 It is clear that the stationary deviation model does not take into account all situations that may arise in practice. For example, widely used physical RNGs can have deviations that are not described by a stationary process (such deviations can be described by slowly varying changes in probabilities, unexpected changes in characteristics, etc.). 
    The computable tests described are quite complex, but  we 
   show that  a simple-to-implement test based on the Lempel-Ziv   code can detect some deviations outside of stationary processes. In particular, for practical testing of PhRNG, 
    we can recommend using a test based on the LZ code (or some other so-called universal dictionary codes) in test batteries.

\section*{Acknowledgment}
Research  was supported  by  Russian Foundation for Basic Research
(grant no. 18-29-03005).


\begin{thebibliography}{10}

\bibitem{Li:08}
Li, M.; Vit{\'a}nyi, P. {\em An Introduction to {K}olmogorov Complexity and
  Its Applications};
\newblock Springer-Verlag:  Berlin/Heidelberg, Germany, 
2008.

\bibitem{Calude:02}
Calude, C.S. {\em Information and Randomness---An Algorithmic Perspective};
\newblock Springer-Verlag:  Berlin/Heidelberg, Germany, 
2002.


\bibitem{downey2010algorithmic}
Downey, R.G.; Hirschfeldt, D.R. {\em Algorithmic Randomness and
  Complexity};
\newblock Springer Science \& Business Media:  Berlin/Heidelberg, Germany, 
2010.

\bibitem{downey2006calibrating}
Downey, R.; Hirschfeldt, D.R.; Nies, A.; Terwijn, S.A. Calibrating
  randomness. {\em Bull. Symb. Log.} \textbf{2006},~\emph{12},~411--491.


\bibitem{merkle2006kolmogorov}
Merkle, W.; Miller, J.S.; Nies, A.; Reimann, J.; Stephan, F. Kolmogorov--loveland randomness and stochasticity. {\em Ann. Pure Appl. Log.} \textbf{2006},~\emph{138},~183--210.


\bibitem{vv}
V'yugin V.
On Nonstochastic Objects
 {\em Probl. Peredachi Inf.},
1985,  21, 2,  3--9

\bibitem{re}
Reimann, Jan. "Information vs. Dimension: An Algorithmic Perspective." Structure And Randomness In Computability And Set Theory (2020): 111.

\bibitem{ve}
Vereshchagin N. Algorithmic Minimal Sufficient Statistics: a New Approach. Theory of Computing Systems. P. 1-19, 2015.


\bibitem{le1}
L'Ecuyer, P. History of uniform random number generation. In Proceedings of the WSC
  2017-Winter Simulation Conference, Las Vegas, NV, USA, 3--6 December 2017.
  

\bibitem{le2}
P.~L'Ecuyer and R.~Simard, ``TestU01: AC library for empirical testing of random number generators,'' {\em ACM Transactions on Mathematical Software}, vol. 33, no. 4, p.22, 2007.



\bibitem{NIST}
A.~Rukhin, J.~Soto, J.~Nechvatal, M.~Smid, E.~Barker, S.~Leigh, M.~Levenson,
  M.~Vangel, D.~Banks, A.~Heckert, J.~Dray, and S.~Vo, {\em A Statistical Test
  Suite for Random and Pseudorandom Number Generators for Cryptographic
  Applications}.
\newblock National Institute of Standards and Technology, 2010.

\bibitem{j4} D. Hurley-Smith, C. Patsakis and J. Hernandez-Castro, On the unbearable lightness of FIPS 140-2 randomness tests,  IEEE Transactions on Information Forensics and Security, doi: 10.1109/TIFS.2020.2988505.



\bibitem{lem}
Lehmann E.L. Statistical Hypotheses. Breakthroughs in Statistics: Foundations and Basic Theory. 2012 Dec 6:67.


\bibitem{ks}
M. Kendall,  A. Stuart,   {\em The advanced theory of statistics; Vol.2:
  Inference and Relationship};
\newblock Hafner Publishing Company: New York, NY, USA, 
 1961.


\bibitem{dr} 
 Ryabko, D.
Asymptotic Nonparametric Statistical Analysis of Stationary Time Series.
Springer, 2019.



\bibitem{R2}
B. Ryabko, J.  Astola, M.  Malyutov,
 \emph{  Compression-Based Methods of Statistical Analysis and Prediction of Time Series,}
  Springer International Publishing Switzerland, 2016.


\bibitem{r-isita}
Ryabko B.  Statistical Testing of Randomness, ISITA2020 — October 24–-27, 2020,  — Kapolei,  USA.
(see also Ryabko B. On asymptotically optimal tests for random number generators. 
https://arxiv.org/abs/1912.06542 )


\bibitem{to}
Tadaki, K. "Partial randomness and dimension of recursively enumerable reals." In International Symposium on Mathematical Foundations of Computer Science, pp. 687-699. Springer, Berlin, Heidelberg, 2009.



\bibitem{rs}
Reimann J, Stephan F. Hierarchies of randomness tests. In: Mathematical logic in Asia 2006 (pp. 215-232).

\bibitem{co}
T.~M. Cover and J.~A. Thomas, \emph{Elements of information theory}.\hskip 1em
  plus 0.5em minus 0.4em\relax New York, NY, USA: Wiley-Interscience, 2006.




\bibitem{rm}
Ryabko, B. Ya, and V. A. Monarev. "Using information theory approach to randomness testing." Journal of Statistical Planning and Inference 133, no. 1 (2005): 95-110.
 
 
\bibitem{re2}
R. G. Downey and J. Reimann.  Algorithmic randomness. Scholarpedia, 2(10), 2007. 
 
 
\bibitem{gram} Charikar, M., Lehman, E., Liu, D., Panigrahy, R., Prabhakaran, M., Rasala, A., Sahai, A. and Shelat, A., 2002, May. Approximating the smallest grammar: Kolmogorov complexity in natural models. In Proceedings of the thiry-fourth annual ACM symposium on Theory of computing (pp. 792-801).




  


\end{thebibliography}

\end{document}